\documentclass[12pt]{article}
\usepackage[utf8]{inputenc}
\usepackage{amsmath}
\usepackage{amsfonts}
\usepackage{amssymb}
\usepackage{amsthm}
\usepackage{relsize}
\usepackage{array}
\usepackage{multirow}
\usepackage{tabu}
\usepackage{soul}
\usepackage{graphicx}
\usepackage{epigraph}
\usepackage{float}
\usepackage[english]{babel}
\usepackage[usenames, dvipsnames]{color}
\usepackage{fancyhdr}
\usepackage[Sonny]{fncychap}
\usepackage{tikz-cd}
\usepackage{tkz-euclide}
\usepackage{sseq}
\usepackage{newfloat}
\usepackage{amsmath,mathtools}
\usepackage{pgfplots}
\usepackage{geometry}
 \usepackage{sectsty}
\sectionfont{\fontsize{14}{15}\selectfont}

\usepackage[all]{xy}
\DeclareFloatingEnvironment[fileext=dia]{diagram}

\theoremstyle{definition}

\newtheorem{remark}{Remark}

\usepackage[most]{tcolorbox}

\makeatletter
\newtcolorbox{note1}[1][]{%
	breakable,
	enhanced jigsaw, 
	sharp corners, 
	boxrule=0pt, 
	attach title to upper, 
	fontupper=\linespread{1.1}\fontfamily{qpl}\selectfont,
	fontlower=\linespread{1.1}\fontfamily{qpl}\selectfont, 
	left=0pt,
	right=0pt,
	top=0pt,
	bottom=0pt,
	boxsep=3pt,
	colback=green!3!white,
	frame hidden,
	before skip=10pt plus 2pt,after skip=10pt plus 2pt,
	#1
}
\makeatother

\usepackage{tabularx}
\usepackage{imakeidx}
\usepackage{hyperref}
\hypersetup{colorlinks={true},linkcolor={blue},citecolor={blue},urlcolor={blue}, linktoc=all}  
 
 \usepackage[capitalize, nameinlink]{cleveref}
 \crefdefaultlabelformat{#2\textbf{#1}#3}
 \crefname{figure}{Figure}{Figures} 
 \Crefname{figure}{Figure}{Figures}
 \crefname{table}{Table}{Tables}
 \Crefname{table}{Table}{Tables}
 \crefname{section}{\S\hspace{-1mm}}{\S\hspace{-1mm}}
 \Crefname{section}{\S\hspace{-1mm}}{\S\hspace{-1mm}}
 \crefname{equation}{}{}
 \Crefname{equation}{}{}
 \crefname{example}{Geometric Pattern}{Geometric Patterns} 
 \Crefname{example}{Geometric Pattern}{Geometric Patterns}

 \usepackage{footnote}

\begin{document}

\title{\textbf{The Elamite Formula for  The   Area of   a Regular Heptagon}}

\author{Nasser Heydari and Kazuo Muroi}

\maketitle

\begin{abstract}
In this article, we   study   the inscription on the  reverse of Susa Mathematical Text No.\,2, a clay tablet held in the collection of the Louvre Museum and thought to date from between 1894--1595 BC. We focus on the formula given in this text for the  approximate area of a regular heptagon. We  give a geometric explanation for the formula and show that this approximation is more accurate than other  contemporaneous formulas in Babylonian mathematics and even that of  Greek mathematician Heron who proved it almost 1800 years later. We also consider the possible ways the Susa scribes might have applied this formula  to construct the regular heptagon for inscription on a clay tablet.
\end{abstract}

\section{Introduction}
This tablet is one of  26 clay  tablets excavated from Susa by French archaeologists in 1933. The texts of all the Susa mathematical tablets (henceforth \textbf{SMT}) along with their interpretations were first published in 1961 (see \cite{BR61}). Since then, many scholars have   reexamined these mathematical texts--including the second author of this article--and the original interpretations have been and continue to be subject to revision.\footnote{Some texts of  \textbf{SMT} have been reexamined in \cite{Fri07-1,Fri07-2,FA16,Hyp93,Hyp99,Hyp02,Hyp02,Hyp17,Mur92-1,Mur92-2,Mur92-3,Mur94-1,Mur94-2,Mur00-1,Mur01-1,Mur01-2,Mur03-2,Mur13,Mur16,Mur}. The authors of this paper are working on a book about Elamite mathematics in which all texts of \textbf{SMT} are reexamined. In this project, we are using the newly photographed and  high-resolution photos of all \textbf{SMT}  which we acquired from the Louvre Museum.}

 The structure of the tablet is as follows\footnote{The reader can see the new  photos of this tablet on the website of the Louvre's collection. Please see \url{https://collections.louvre.fr/ark:/53355/cl010185652} for obverse  and \url{https://collections.louvre.fr/ark:/53355/cl010185652}, for reverse.}: there is a  regular hexagon    on the obverse of this tablet with some numbers, and  on the  reverse there is a regular heptagon  with numbers and  a formula  for its area. The formula gives an instruction to find the area constant of a regular heptagon. Unfortunately, this formula has  been given little attention by scholars and many of them have only mentioned it as providing an approximation for the area of a regular heptagon. However, we   believe this formula deserves deeper consideration because it uses  geometric ideas to give  a very good approximation and accuracy led us  to write this article   to discuss fully    the  mathematical aspects of the formula.

 It is also noteworthy that  the dimensions of the  two polygons inscribed on  both sides of the tablet under consideration are   accurate. Also,    a part of  the circumscribed circle  is still visible on both sides of the tablet suggesting that  the Susa scribe    used  both a   compass and a straightedge to draw the   figure. 

\begin{remark}
By convention,  we have used the sexagesimal numeral system to write numbers. In this system, the comma ``,'' separates the double digits  and the   semicolon  ``$;$''    separates the non-negative and negative powers of 60. So, a number like $12,23,5;13,45,9$ in the sexagesimal numeral system becomes   
$$12\times 60^2+23\times 60^1+5\times 60^0+13\times 60^{-1}+45\times 60^{-2}+9\times 60^{-3}$$
 in the decimal numeral system.     
\end{remark}
 
\section{Areas of polygons in Babylonian and Elamite mathematics}
Let $n\geq 3$ be a natural number. A  \textit{polygon} with $n$ sides or an $n$\textit{-gon} is a plane figure formed by a chain of $n$  line segments $v_1v_2$, $v_2v_3$, $\cdots$, $v_{n-1}v_{n}$ and $v_nv_1$ connecting $n$ points $v_1,v_2,\cdots, v_n $ and enclosing a bounded region of the plane. In such a case, one may call the bounded plane region, the boundary chain, or the two together   as  the polygon. The line segments $v_{i}v_{i+1}$ are the \textit{sides} and the points  $v_i$ are the \textit{vertices} of the polygon.   The \textit{centroid} of a polygon is the arithmetic mean position of all the points in the figure.

\begin{figure}[H]
	\centering
	\includegraphics[scale=1]{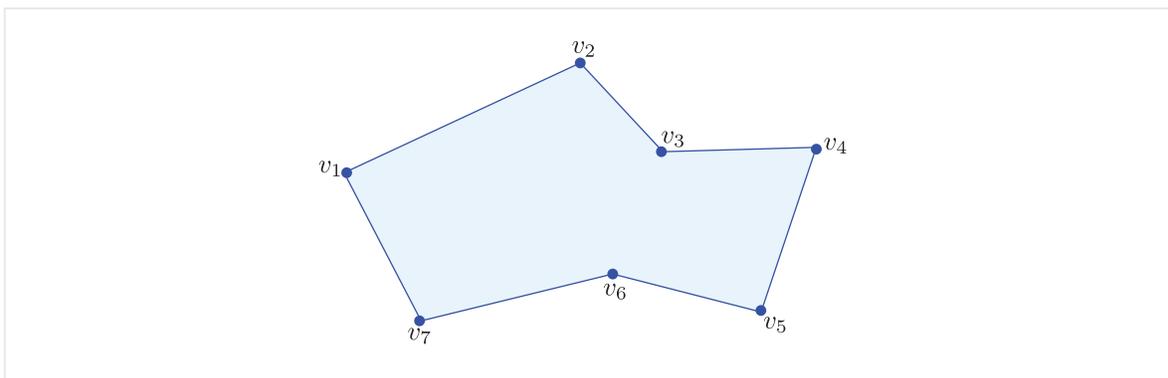}
	\caption{A polygon a with 7 sides and 7 vertices}
	\label{figure1b}
\end{figure}  

A  \textit{regular   polygon}  is an $n$-gon which is equiangular (all of its inner angles have the same measure $ \frac{360^{\circ}}{n}$) and equilateral (all of its sides are of the same length $a$). We use the notation $\Gamma_n(a) $  for a regular polygon with $n$ sides of equal length $a$.   Any regular  polygon can be inscribed in a unique circle whose center coincides with the centroid  of the regular polygon and    passes through all vertices of the polygon. This unique circle with radius $r$ is called the \textit{circumscribed circle} of the regular polygon. Any straight line connecting a vertex of the regular polygon to its centroid is a \textit{circumradius} of the regular polygon. The height from the center to any side is called  the  \textit{apothem} of the polygon.

\begin{figure}[H]
	\centering
	\includegraphics[scale=1]{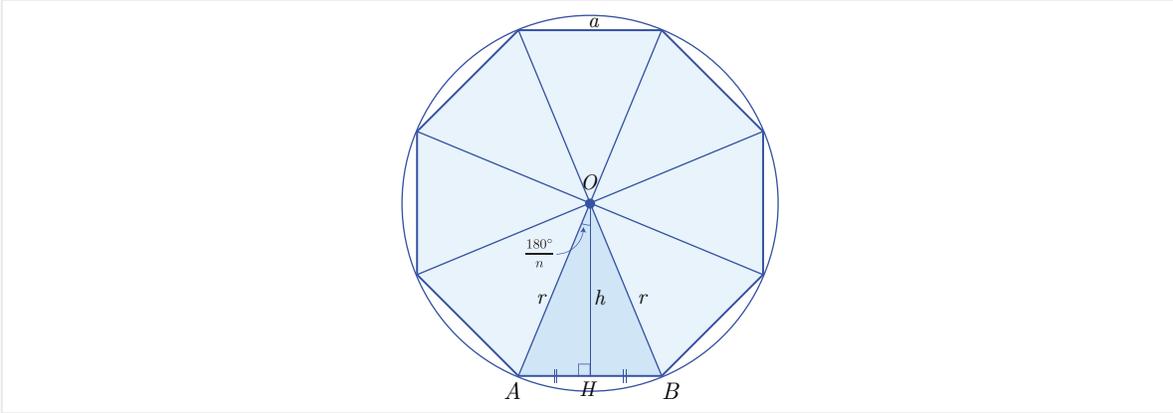}
	\caption{A regular polygon and its circumscribed circle}
	\label{figure1}
\end{figure}  

 Elamite and Babylonian  scribes were seemingly interested in computing the area  of   regular $n$-gons. This task seemed to be  easy for   them when      $n=3$ or $n=4$. In fact, for $n=4$  they    had a square of side $a$ whose area  was  simply computed as   $S_{\Gamma_{4}(a)}= a^2$ while for $n=3$, the polygon was an equilateral triangle with side $a$ whose area could have been computed by the usual formula ``half of base times height''\footnote{It is a well-known fact that ancient scribes used the known standard rules for computing the areas of basic figures: ``length times width'' for rectangles, ``square of side'' for squares, ``half of base times height'' for   triangles, ``product of diagonals'' for rhombuses, and ``half of height times the sum of bases'' for trapezoids (see \cite{Fri07-1}, for example).}.  
 
 \begin{figure}[H]
 	\centering
 	\includegraphics[scale=1]{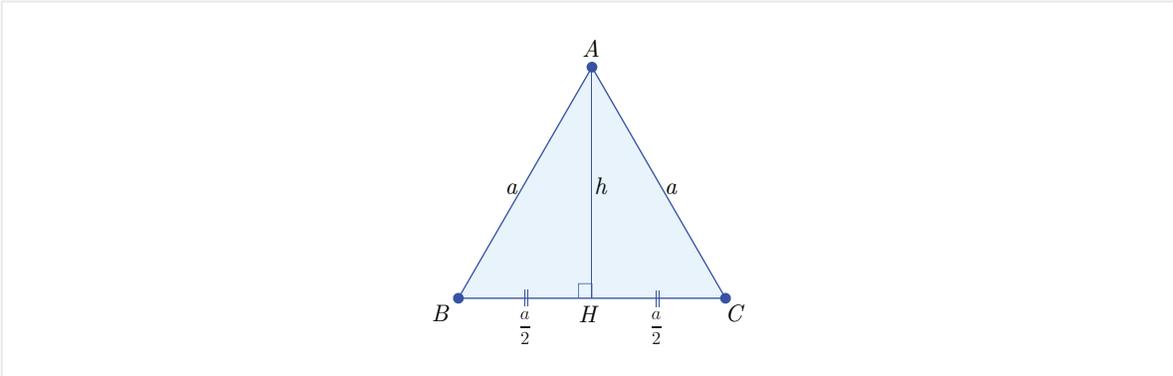}
 	\caption{An equilateral triangle}
 	\label{figure1a}
 \end{figure}  

 It is believed that Babylonian scribes have used the Pythagorean theorem and the approximation  $ \sqrt{3}\approx \frac{7}{4}$   for an equilateral triangle $\Gamma_{3}(a) $ with side $a$ to compute its approximate area. Because in that case, the height is obtained by the Pythagorean theorem as  $h^2=  a^2-\frac{a^2}{4}$ (see \cref{figure1a}). So,  $h =\frac{\sqrt{3}}{2}a\approx \frac{7}{8}a$ and thus the area $\frac{ah}{2} $ becomes
 \begin{equation}\label{equ-SMT2-adaa}
	S_{\Gamma_{3}(a)}\approx\dfrac{7}{16}a^2.
\end{equation}

 Although the area formulas for equilateral triangles and squares could have been easily obtained by scribes, the trouble would have started when they would increased the number of sides , i.e., for the cases $n\geq 5$. As has been discussed by many scholars, the ancient scribes might have   utilized the  known standard method    to compute the areas of polygons, namely  using a circumscribed circle. Recall that to find the area of a  regular $n$-gon with side $a$ and  circumradius $r$, we divide it into $n$ equal  isosceles triangles and try to find the area of this common isosceles triangle whose top angle is $\frac{360^{\circ}}{n} $ (see \cref{figure1}). The general formula  for the exact area of a regular polygon   $ \Gamma_n(a)$ is easily obtained as
 \begin{equation}\label{equ-SMT2-c}
 	S_{\Gamma_n(a)}=  \frac{na^2}{4}  \times \cot\left(\frac{180^{\circ}}{n}\right).
 \end{equation}

\section{\textbf{SMT No.\,2} and Elamite formula}
\begin{figure}[H]
	\centering
	\includegraphics[scale=1]{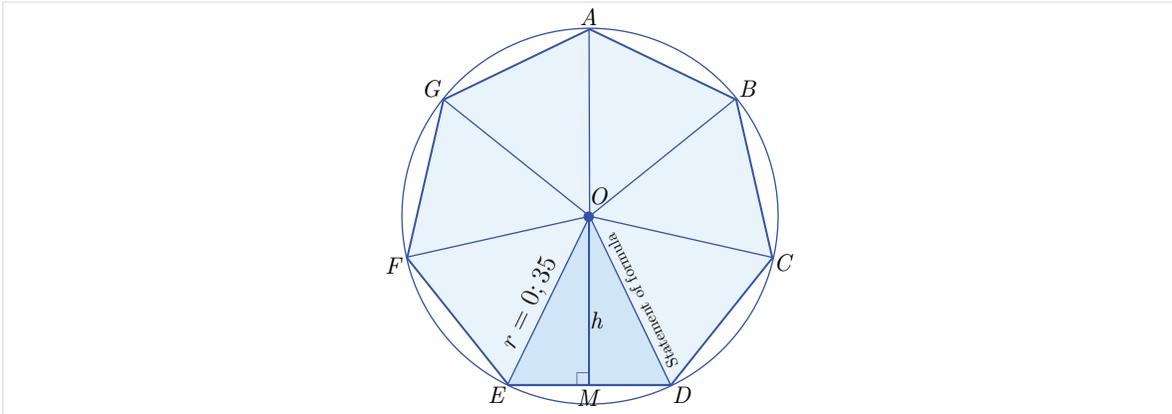}
	\caption{Reconstruction of the  reverse of \textbf{SMT No.\,2}}
	\label{figure2}
\end{figure}

    In \cref{figure2}, we have reconstructed  the   reverse of \textbf{SMT No.\,2}  in which a regular heptagon   $ \Gamma_7(a): ABCDEFG$  is inscribed in a circle  with radius $r=0;35 $. If we connect the center of the circle  to each vertex, we    get  seven equal isosceles triangles  whose  bases  appears to  have    the same value $0;30$. In fact,  by approximating the circumference of the circle    $c_{\Gamma}=2\pi r $  with that of the regular heptagon    $ c_{\Gamma_7}=7a$ and using the Babylonian approximation $\pi\approx 3 $, we have $7a\approx 6\times(0;35)$  or $a=  0;30$.

Unlike the case of regular hexagon on the obverse of this tablet   in which the side $a$ is given, here the circumradius $r$ is known. In fact, the number $0;35$ is clearly recognizable as
\begin{center}
	 {\fontfamily{qpl}\selectfont 35 u\v{s}} ``$0;35$ is the length'' 
\end{center}
  over one of the legs of the isosceles triangle $\triangle OED $, namely $\overline{OE}=0;35 $.

Le us  use this  data and the Babylonian method  to compute the area of the isosceles triangle $\triangle OED $ and the regular heptagon $\Gamma_7(0;30)$. To do so, consider the isosceles triangle $\triangle OED $ as shown in \cref{figure2}. Since the height   $OM$    bisects the base $ED$,  $\overline{EM}=\overline{DM}=\frac{0;30}{2}=0;15 $. By using the Pythagorean theorem   in the right triangle  $\triangle OME $, we have
\begin{align*}
	h&= \sqrt{\overline{OE}^2-\overline{EM}^2}\\
	&\approx\sqrt{(0;35)^2-\dfrac{(0;30)^2}{4}}\\
	&=  \sqrt{\frac{1}{4} \Big[ 4\times (0;20,25) -  0;15  \Big]}\\
	&= \frac{1}{2}\sqrt{1;21,40- 0;15}  \\
	&=\frac{1}{2} \sqrt{1;6,40}.
\end{align*}
Hence, the height $h$ is 
\begin{equation}\label{equ-SMT2-ae}
	h\approx \frac{1}{2} \sqrt{1;6,40},
\end{equation}

To get the approximate value for $h$,  we need to apply  the approximate formula  $\sqrt{1+x}\approx 1+\frac{x}{2} $   to \cref{equ-SMT2-ae} which is the Babylonian method.\footnote{It is a fact that Babylonians used the formula $\sqrt{a^2\pm b}\approx a\pm \frac{b}{2a}$ for approximating irrational square roots like $\sqrt{2}, \sqrt{3},\sqrt{5}$ and so on. See \cite{FR98}, for a discussion on this topic.} Thus, we can continue as follows: 
\begin{align*}
	h&\approx  \frac{1}{2} \sqrt{1;6,40}\\
	& = \frac{1}{2}   \sqrt{1+0;6,40}\\
	& \approx \frac{1}{2} \times \left(1+\dfrac{0;6,40}{2}\right) \\
	& = \frac{1}{2} \times (1;3,20) \\
	&=0;31,40.
\end{align*} 
So, we get the following approximate value of the height $h$  as
\begin{equation}\label{equ-SMT2-a}
	h\approx 0;31,40.
\end{equation}
This easily implies the area of the isosceles triangle $ \triangle OED $ as follows:
\begin{align*}
	S_{\triangle OED}&=\frac{1}{2}\times \overline{ED}\times h \\
	& \approx \frac{1}{2} \times (0;30)\times (0;31,40) \\
	& = (0;15) \times (0;31,40) \\
	&=0;7,55.
\end{align*}
That is, 
\begin{equation}\label{equ-SMT2-b}
	S_{\triangle OED}\approx 0;7,55.
\end{equation}
Finally, we can compute the approximate area of the regular heptagon $\Gamma_7(0;30) $ by using      \cref{equ-SMT2-b}:
\begin{align*}
	S_{\Gamma_7(0;30)}  &= 7\times S_{\triangle OED}\\
	&\approx  7\times  (0;7,55)\\
	&  = 0;55,25.
\end{align*}    
Therefore, we get    
\begin{equation}\label{equ-SMT2-c}
	S_{\Gamma_7(0;30)}\approx 0;55,25.
\end{equation} 

An  approximate formula for the area of a general regular heptagon $\Gamma_7(a) $ of side $a$ can  be  obtained from  \cref{equ-SMT2-c}. In fact,  since $a=2a\times (0;30)$,  the area of the regular heptagon $ \Gamma_7(a)$ can be computed as follows:
\begin{align*}
	S_{\Gamma_7(a)}&= (2a)^2\times S_{\Gamma_7(0;30)} \\
	&  \approx 4\times (0;55,25)\times a^2\\
	&  =(3;41,40)\times a^2.
\end{align*}
If we ignore the second sexagesimal numeral 40, we get  the following Babylonian approximate formula
\begin{equation}\label{equ-SMT2-e}
	S^B_{\Gamma_7(a)}\approx (3;41)\times a^2.
\end{equation}
If $a=1$, the approximate area becomes $ 3;41$ which is  the very geometric  coefficient $3;41$   listed in \textbf{SMT No.\,3}  as the constant of a regular heptagon. This confirms that \cref{equ-SMT2-e} would have been used to approximate the area of a regular heptagon.

Although we do not know whether the previous numbers in   \cref{equ-SMT2-a}, \cref{equ-SMT2-b} and \cref{equ-SMT2-c}  have been written down on the missing part of the tablet, we can read a formula for the area  of a regular heptagon  below the missing  part:

\begin{note1} 	
	\underline{Reverse}
	
	(L1) [s]ag-7 \textit{a-na} 4 \textit{te-}[\textit{\c{s}i}]\textit{-ip-ma}\\
	(L2)  \textit{\v{s}\'{i}-in-\v{s}\'{e}-ra-ti}\\
	(L3) \textit{ta-na-as-s\`{a}-ah-ma} a-\v{s}\`{a}
\end{note1}

	\noindent
	\underline{Translation:}\\
	``A (regular) heptagon. You multiply (the square of a side) by 4 and you subtract one-twelfth (of the result from the result itself), and (you see) the area.'' \\

In other words, to find the area     of a regular heptagon  $\Gamma_7(a)$ with side of length $a$ one needs to multiply $a$ by 4 and then subtract the twelfth of $4a$ from itself. If we translate this  instruction   in mathematical language and perform the related  calculations, we get: 

\begin{align*}
	S^E_{\Gamma_7(a)}  &= 4\times a^2 - \frac{1}{12}\times \left(4\times a^2\right)\\
	&=   \left(4-\frac{4}{12}\right) \times a^2 \\
	&=   \left(\frac{48-4}{12}\right) \times a^2 \\
	&=   \frac{44}{12}    a^2 \\
	&=    \frac{11}{3}  a^2 \\
	&=  (3;40)\times a^2.
\end{align*} 

\noindent
This gives us another   formula for the area of a regular heptagon as follows
\begin{equation}\label{equ-SMT2-d}
	S^E_{\Gamma_7(a)}= (3;40)\times a^2, 
\end{equation} 
which from now on we call the \textit{Elamite formula}.

Although the Elamite formula \cref{equ-SMT2-d} seems to be  only another approximate  formula for the area of a regular heptagon,  it provides us with a good approximate value which is more accurate than  not only  the Babylonian  formula g  \cref{equ-SMT2-e} but also    the Greek formula   (see \cite{Hea21}, Page 328) due to  Heron   of Alexandria  (circa 50 AD):
\begin{equation}\label{equ-SMT2-ea}
	S^H_{\Gamma_7(a)}= \frac{43}{12} a^2 = (3;35)\times a^2.
\end{equation}  

To see the accuracy of these approximate formulas, we compute their error percentages.  First, note  that  by \cref{equ-SMT2-c}  the accurate area  of the regular heptagon  $\Gamma_7(a)$ to four sexagesimal places is
\begin{align*}
	S^A_{\Gamma_7(a)}&=  \frac{7}{4}\times \cot\left(\frac{180^{\circ}}{7}\right)\times a^2\\
	&= (1;45) \times  (2;4,35,28,37,17,\cdots)\times a^2\\
	&=    (3;38,2,5,5,16,\cdots)\times a^2\\
	&\approx   (3;38,2,5,5)\times a^2,
\end{align*}
then  observe that

\begin{align*}
	e^H_7&= \left|\frac{S^A_{\Gamma_7(a)}-S^H_{\Gamma_7(a)}}{S^A_{\Gamma_7(a)}}\right| \times 100\% \\ &\approx\frac{(3;38,2,5,5-3;35)\times a^2}{(3;38,2,5,5)\times a^2} \times 100\% \\
	&= \frac{0;3,2,5,5}{3;38,2,5,5} \times 100\% \\
	&\approx (0;0,50,6,25) \times 100\%\\
	&\approx   1.39\%, 
\end{align*} 
and
\begin{align*}
	e^B_7&= \left|\frac{S^A_{\Gamma_7(a)}-S^B_{\Gamma_7(a)}}{S^A_{\Gamma_7(a)}}\right| \times 100\% \\ &\approx\frac{(3;41-3;38,2,5,5)\times a^2}{(3;38,2,5,5)\times a^2} \times 100\%  \\
	&= \frac{0;2,57,54,55}{3;38,2,5,5} \times 100\% \\
	&\approx (0;0,48,57,35) \times 100\%\\
	&\approx  1.36\%, 
\end{align*} 
as well as
\begin{align*}
	e^E_7&= \left|\frac{S^A_{\Gamma_7(a)}-S^E_{\Gamma_7(a)}}{S^A_{\Gamma_7(a)}}\right| \times 100\% \\ 
	&\approx\frac{(3;40-3;38,2,5,5)\times a^2}{(3;38,2,5,5)\times a^2} \times 100\% \\
	&= \frac{0;1,57,54,55}{3;38,2,5,5} \times 100\%  \\
	&\approx (0;0,32,26,55) \times 100\%\\
	&\approx 0.9\%. 
\end{align*} 
It is clear that among three values $ e^H_7, e^B_7$ and  $e^E_7$, the last one is the smallest error percentage  proving that formula \cref{equ-SMT2-d} produces the most  accurate value for the area of a regular heptagon. In fact, the error is less than one percent which is remarkable and shows  the impressive work of the Susa scribes.

\section{Geometric Explanations for Elamite formula} 

Unlike Heron's formula, we do not know how the Susa scribes   arrived at the  formula  \cref{equ-SMT2-d}  inscribed on this tablet. It seems that Heron\footnote{For more details, see \cite{Hea21}.} utilized two approximations $ r\approx \frac{8}{7}a$ and  $ \sqrt{23}\approx \frac{43}{9}$   for deriving formula \cref{equ-SMT2-ea}:
\begin{align*}
S^H_{\Gamma_7(a)}&\approx 7 \times \frac{a}{2}\times \sqrt{\frac{64a^2}{49}-\frac{a^2}{4}} \\
&=\frac{7a^2}{2}\sqrt{\frac{9\times 23}{4\times 49}} \\
&=\frac{3\sqrt{23}}{4}a^2 \\
&  \approx \frac{3}{4} \times   \frac{43}{9}\times a^2\\
&= \frac{43}{12}a^2.	
\end{align*}

It is interesting that Heron has used the approximation $ a\approx \frac{7}{8}r$ to construct a regular heptagon of side $a$  inscribed in a circle of radius $r$. In other words, he has used the apothem of the inscribed regular hexagon as the side of the regular heptagon (see \cref{figure11}). As we will see in \cref{sec-construction}, this choice turns out to be a pretty good approximation.  

Although Susa scribes have not explained how they have found their formula,  one may use    geometric  explanations   to obtain the Elamite formula  for  the area of a regular heptagon  with side $a$.  We try to give our explanation. 

According to the statement of the formula,   one first needs to consider the area $4a^2 $ and then  subtract   $ \frac{1}{12}(4a^2)$ from  it in order to find the approximate formula of a regular heptagon. This suggests that we first need to consider a geometric figure with area $4a^2$ which can be easily divided into 12 equal parts. While there are many different geometric shapes  whose area is $4a^2$, the second condition can lead us to choose one of the  following basic figures, i.e., \\
1) a square with side $2a$; or\\
2) a rectangle with sides $a$ and $4a$. 

\begin{figure}[H]
	\centering
	\includegraphics[scale=1]{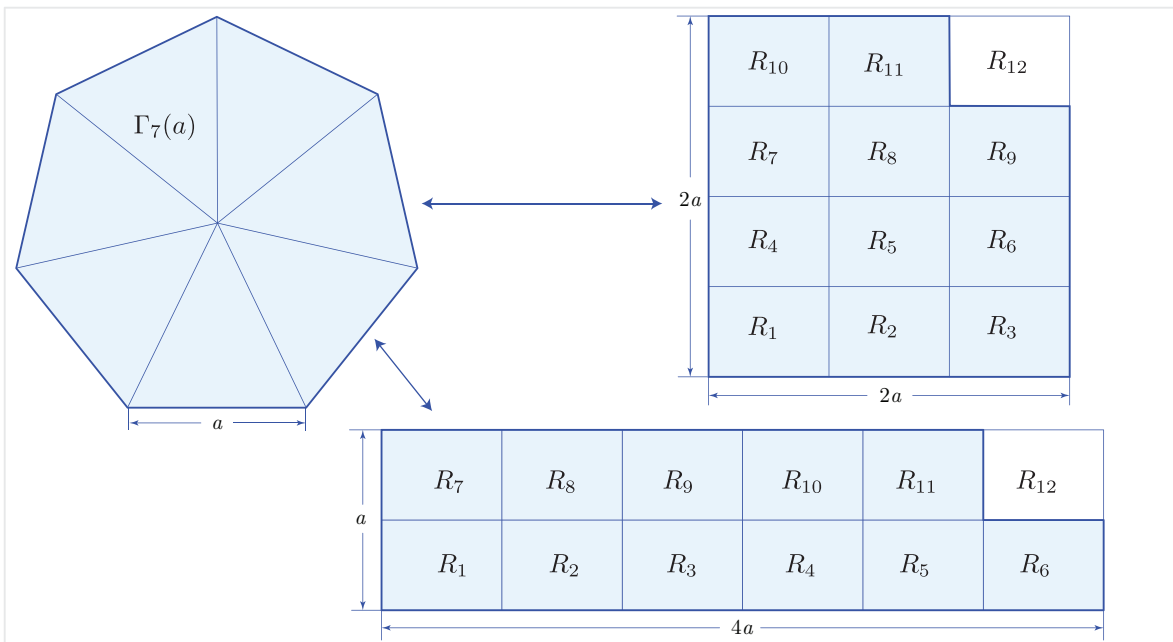}
	\caption{Two approximations for a regular heptagon}
	\label{figure3}
\end{figure}

The main advantage of these choices is that we can easily divide them into equal parts by using vertical and horizontal lines. As is shown in \cref{figure3}, a square and a rectangle are divided in 12 equal small rectangles $R_1, R_2,\ldots, R_{12}$ by horizontal and vertical lines. In each case, if    one-twelfth of the shape is removed,  the remaining part is claimed to be an approximation for the area of the regular heptagon with side $a$.    We have removed the last rectangle $R_{12}$   in each figure.

We use a geometrical approach, which we would like to call \textit{cut and paste}, to verify the Elamite formula. This method is very useful when we want to show that two figures have the same (or almost the same) area. We usually have  an original figure and a goal figure whose area approximates that of the original figure. The main idea of this method is   to   cut  the original figure into smaller pieces and then to paste them together in a  way  so as to  (almost) fill up the goal figure. The most obvious example for this method might be an isosceles triangle and a rectangle. If we cut the isosceles triangle along its height, we get two equal right triangles. By attaching these two right triangles along their hypotenuses, we obtain a rectangle with the same area as the original isosceles triangle. In the current problem, the original figure is a regular heptagon  and the goal figure is either a square  without a corner or a rectangle  without a corner (see \cref{figure3}).

\begin{figure}[H]
	\centering
	\includegraphics[scale=1]{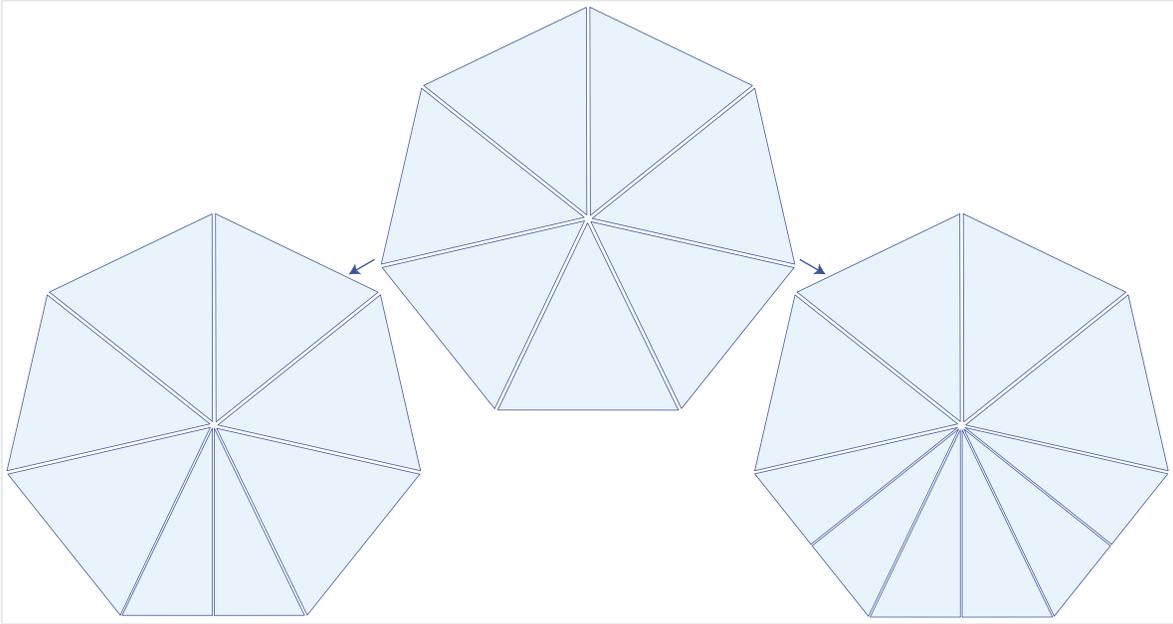}
	\caption{Cutting off a regular heptagon in two ways}
	\label{figure4}
\end{figure}

Start with   a regular heptagon  $\Gamma_7(a)$    and  connect its center to each vertex by straight lines. Then slice off the regular heptagon  $\Gamma_7(a)$  along those lines to get seven equal isosceles triangles  whose bases are of  length  $a$. Now, cut off two or four of these isosceles triangles  along their heights to get four or eight  equal right triangles. So  we  obtain 9 or 11   triangles   as shown in      \cref{figure4}.

\begin{figure}[H]
	\centering
	\includegraphics[scale=1]{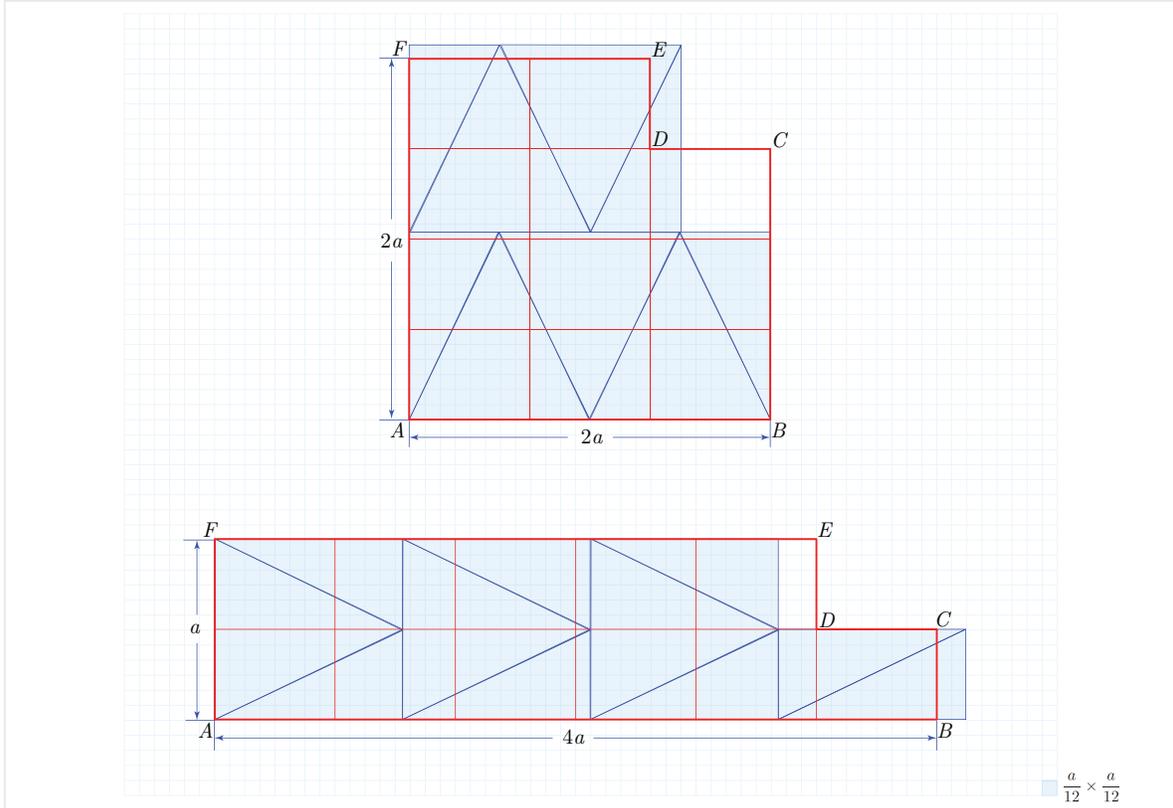}
	\caption{Square and rectangular layouts of  triangular pieces}
	\label{figure5}
\end{figure}

Now,   consider a grid pattern of squares whose sides are of length  $\frac{a}{12}. $\footnote{In fact, any multiple of 12 works and the bigger the number, the better the approximation.} We   consider our   figures in this grid pattern. For each figure, we try to lay out the obtained triangles  on the grid so as to fill up the figure      (see \cref{figure5}).\footnote{In this picture we have scaled the dimensions to help the reader clearly see what is happening.} For the square without corner, there are five isosceles triangles and four right triangles, while  for the rectangle without corners we have used three isosceles triangles and eight right triangles.  Note that in each case a part of the colored area goes out of the red hexagon $ABCDEF$  and another part of it is not covered by colored area.

Finally, we cut the small colored  squares outside the   hexagon $ABCDEF$ and set  them out  in the blank region inside it. In fact,  we can easily count the colored and blank small   squares and rectangles in each case.  \\
\textbf{Case 1. Square layout}\\
According to the upper layout in \cref{figure5}, we count as follows: \\ 
i) complete small colored squares: 12\\
ii) almost complete small colored squares: 18\\
iii) complete blank squares: 30\\
iv) almost blank half squares: 6 (= 3 almost complete blank squares).\\
So there are totally 30 colored squares and 33 blank ones in this case.     After arranging the colored ones inside the blank region, there remains three blank small squares. 

\noindent
\textbf{Case 2. Rectangular layout}\\
According to the lower layout in \cref{figure5}, we count as follows: \\ 
i)  almost complete small colored squares: 12\\
ii) complete blank squares: 12\\
iii) almost blank half squares: 6 (= 3 almost complete blank squares).\\
So there are totally 12 colored squares and 15 blank ones in this case.     After arranging the colored ones inside the blank region, there remains three blank small squares.

\begin{figure}[H]
	\centering
	\includegraphics[scale=1]{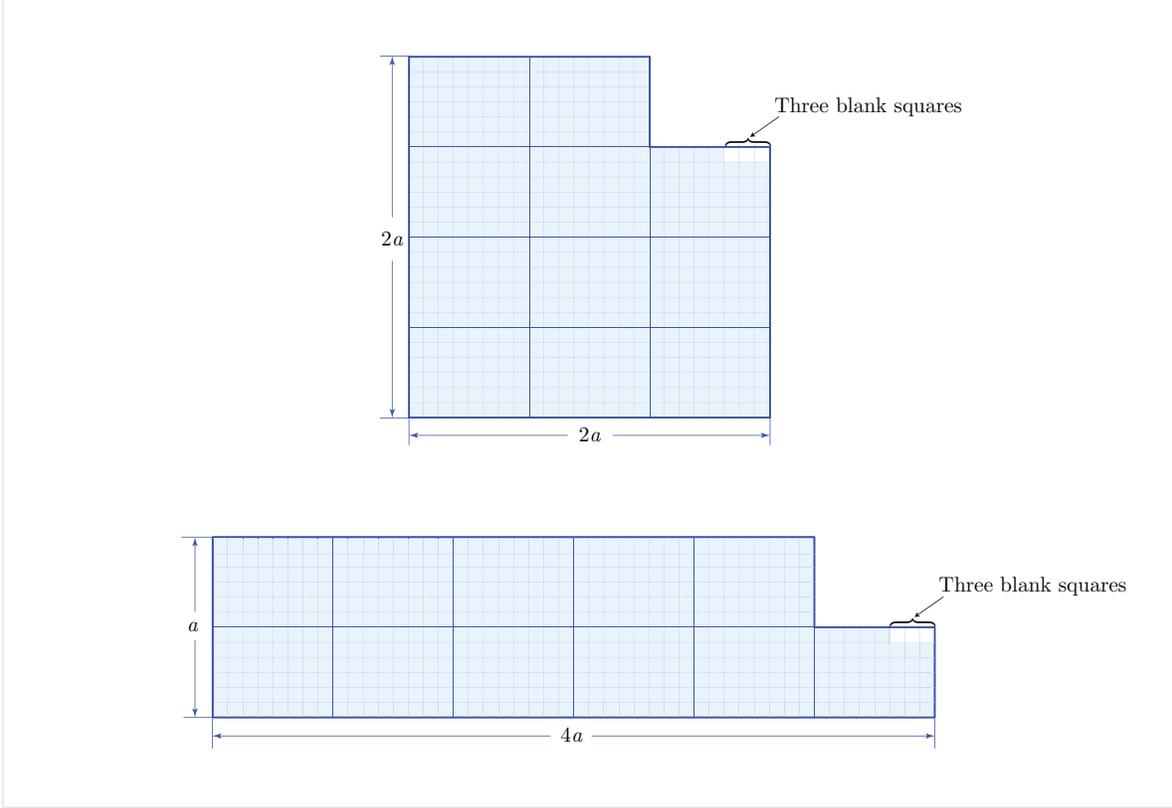}
	\caption{Two layouts of a regular heptagon}
	\label{figure6}
\end{figure}

In both case,  there are three blank small square of side $\frac{a}{12}$ left and  their total areas   is (see \cref{figure6})  
\[ 3\times \left( \frac{a}{12}\right)^2=\frac{3a^2}{144}=\frac{a^2}{48}. \]
Therefore, we have
\begin{align*}
	S_{\Gamma_7(a)} =  S_{ABCDEF} -\frac{1}{48}a^2\approx S_{ABCDEF}.
\end{align*}
Note that the error in this approximation  is 
\[\left| S_{\Gamma_7(a)}-S_{ABCDEF}\right|=\frac{1}{48}a^2\]
and the  error percentage is
\[\frac{\left| S_{\Gamma_7(a)}-S_{ABCDEF}\right|}{S_{ABCDEF}} \times 100\% =\left( \dfrac{\ \ \dfrac{1}{48}a^2\ \ }{\ \ \dfrac{11}{3} a^2\ \ }\right) \times 100\%=\frac{3}{528}\times 100\%\approx 0.57\%,\]
which is  surprisingly small. 

By  comparing \cref{figure5} with \cref{figure6}, the reader may appreciate how the  Susa scribes  may have  used the Elamite formula \cref{equ-SMT2-d} to approximate the area of a regular heptagon. It should be emphasized that we cannot state with certainty   the  reasoning behind the Elamite formula, only that in our view there seems to be considerable   geometric intuition   behind the formula  confirming  that the Susa scribes  must have possessed a great deal of mathematical skill  and experience, which  enabled them  to come up with such a beautiful and accurate formula.

\begin{figure}[H]
	\centering
	\includegraphics[scale=1]{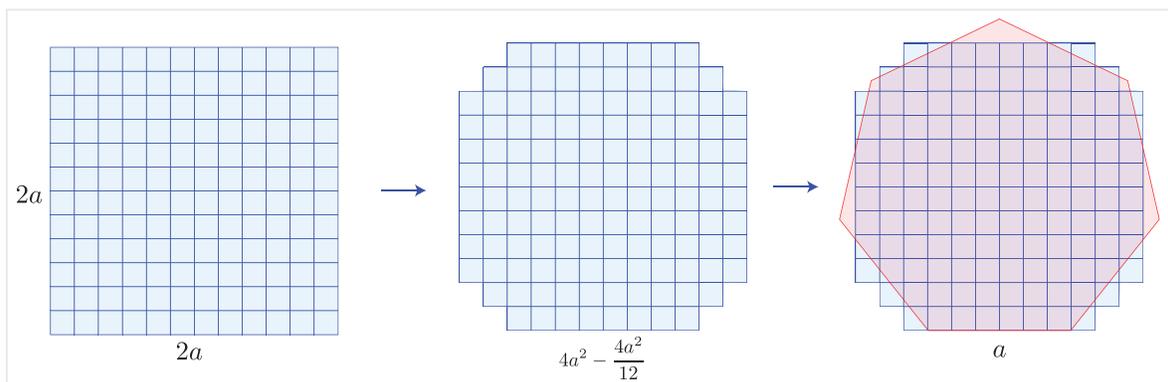}
	\caption{Another visual representation for Elamite formula}
	\label{figure15a}
\end{figure}

\begin{remark}
The reader should note that the geometric explanation we have provided is one of several possible explanations and there are definitely  other ways to justify this formula. For instance,  \cref{figure15a} shows another visual representation for this approximation.
\end{remark}

\section{Significance of Elamite formula}
By way of addition to the history of mathematics, this formula might be of great importance in showing that Elamite scribes    not only had a great deal of mathematical skill  but also were masters of geometry. It seems that geometry    played a major role in the lives of Elamite people, in particular in those of Elamite artisans, as the  wide application of geometric patterns attest  in the Elamite art.\footnote{The reader can consult the   book    ``Art of Elam'' in which there are   many colored plates representing Elamite artifacts  \cite{Alv20}.}  Unfortunately, this social and cultural  characteristic of Elamite civilization  has not been fully appreciated by other scholars who merely consider the Susa mathematical tablets as a part of Babylonian mathematics. The origin of  the Elamite civilization as an independent  entity in southwestern   Iran  dates more or less to (3500--3200 BC)    that of the Sumerians in Mesopotamia. Like other ancient civilizations in the Near East, Elamites had their own culture, customs,  traditions and  writing systems one of which is  recently argued  to be the first  known phonetic writing system.\footnote{The reader can consult the newly published paper  \textit{The Decipherment of Linear Elamite Writing}   \cite{DTKBM} whose authors have deciphered this ancient  inscription for the first time. This writing system writing was used in southern Iran in the late third/early second millennium BC (2300--1880 BC).}  Elamites were   western  neighbors to the Mesopotamian civilizations (Sumer, Akkad, Babylonia, and Assyria) for almost three   millennia and had cultural, social, and  literary  interactions with them. Further study of the \textbf{SMT} may revise the independent contribution of the Elamite civilization to mathematics in the ancient world.

The Elamite formula for the area of a regular heptagon demonstrates   the  skill  of Susa scribes in   computing the areas of geometric figures. However   other scholars  have not appreciated its significance. For example, in \cite{Fri07-1, Fri07-2, PAT92, FR98},   this formula is considered to produce only   another Babylonian constant for the area of a regular heptagon. In \cite{BR61},   only   the transliteration of the formula is given  without any further information.   Most scholars think that the constant $3;40$ is an approximate value of the more familiar constant $3;41$. They do not consider it an independent approximation and  ignore the fact  that it can be derived from a geometric idea. A strong  reason to oppose this view is that the Susa scribe   has given an instruction  on the tablet that produces this constant.     Additionally, there is  the impressive accuracy of this   formula and the geometric beauty which underlies it.  By studying the geometric patterns in the Elamite art alongside  the mathematical ability  of the Susa scribes to compute the areas of complicated figures\footnote{See  the list of geometric constants  given in  \textbf{SMT No.\,3}. In this text, the Susa scribes have provided a list of  constants related to parts or   areas of geometric figures among which there are complex figures such as a  ``concave square '' and ``regular concave hexagon''.}, we   came to the   conclusion that this formula was independently derived by using geometric skills and mathematical experience developed within the Elamite civilization.

\section{Construction of regular polygons}\label{sec-construction}
To construct a regular $n$-gon by compass and straightedge, it is usually very helpful to use a circle, because    in this method we just need  to divide the circumference of the circle into $n$ equal arcs. After doing that, we  connect the adjacent   points to get $n$ equal chords of the circle which are the sides of the required regular $n$-gon.

 \begin{figure}[H]
 	\centering
 	\includegraphics[scale=1]{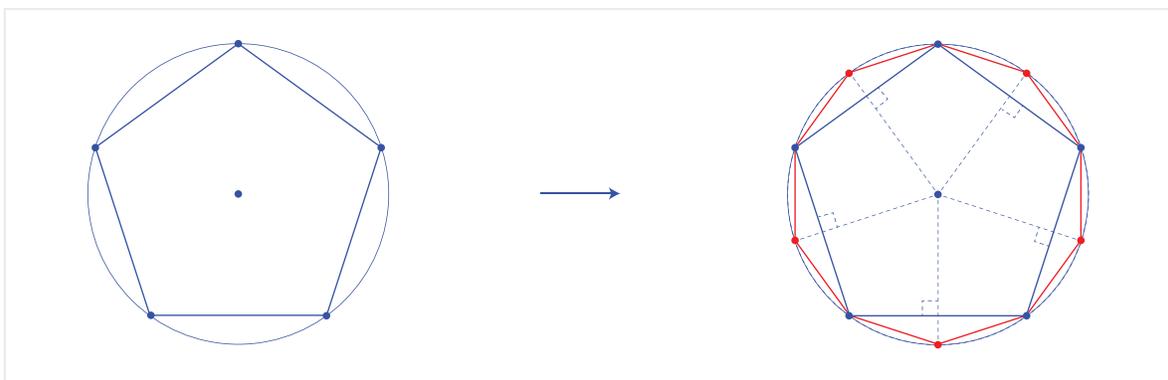}
 	\caption{Construction of a regular $2n$-gon}
 	\label{figure7}
 \end{figure}

We can construct new regular polygons out of given ones.  There are two basic methods to do so. Firstly, when we have a regular $n$-gon, we can easily find the midpoint of each chord by using straight  lines passing through the center and perpendicular to each side of $n$-gon. These new midpoints along with the vertices of the main $n$-gon can be used to construct a regular $2n$-gon (see \cref{figure7}). By repeating this process, one can construct any $2^{m} n$-gon for any natural number $m$.

For the second method, let  $n$ and $m$ be two   coprime natural numbers where $n>m>2$. Also assume that we have constructed a regular $n$-gon and a regular $m$-gon in the same circumscribed circle. We can use these two polygons to construct a regular $nm$-gon as follows:\\
(1)  First construct a regular $n$-gon with vertices $v_{1},v_{2},\cdots,v_{n}$.\\
(2) Construct a regular $m$-gon with vertices $w_{1},w_{2},\cdots,w_{m}$ where $w_{1}=v_{1}$, i.e., they have a common vertex.\\
(3) Repeat step 2  with $v_{2},\cdots,v_{n} $ as the first vertex each time.\\
(4) Mark the $mn$ vertices obtained and connect the adjacent ones to form a regular $mn$-gon.  

\begin{figure}[H]
	\centering
	\includegraphics[scale=1]{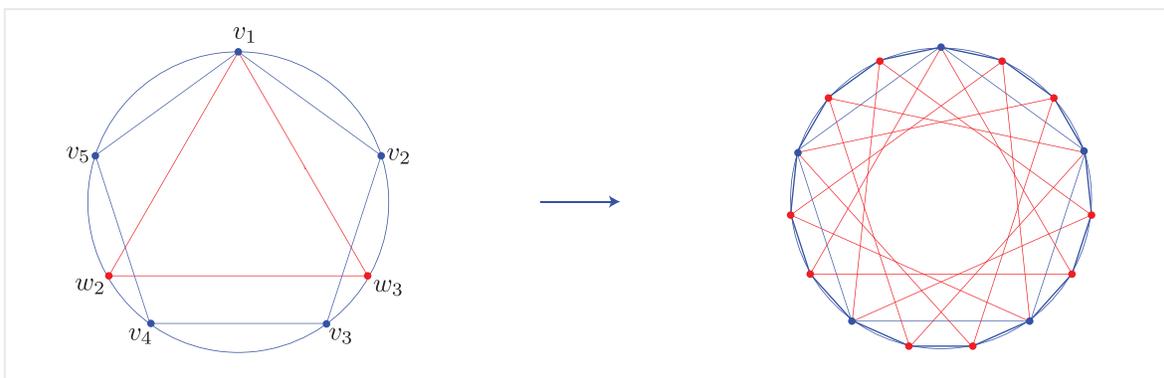}
	\caption{Construction of a regular $mn$-gon}
	\label{figure8}
\end{figure}

The key point in the second method is that the obtained vertices in each step never coincide with the previous ones and we end up with exactly $nm$ points separating equal arcs of the circle (see \cref{figure8}).

Among all the basic regular polygons, the construction of a regular hexagon is probably  the  easiest, because its side and the radius of its circumscribed circle are equal. This unique characteristic enables us to use the same circle to divide the circumference of the circumscribed circle into 6 equal arcs. We can do so as follows: choose an arbitrary point   on a circle of radius $r$. Then put the compass to distance $r$ in the point and draw an arc  intersecting the circumscribed circle in a new point. Repeat this process five more  times to get six  points and  connect  adjacent points to get a regular hexagon.

\begin{figure}[H]
	\centering
	\includegraphics[scale=1]{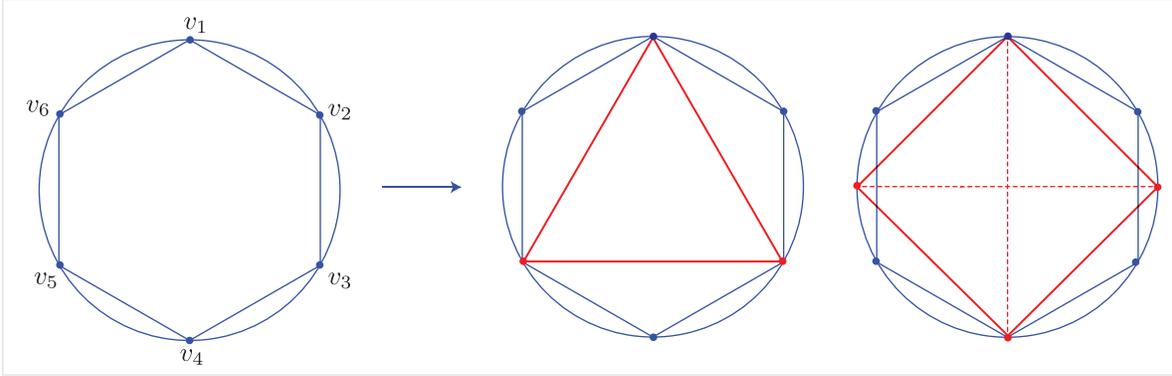}
	\caption{Construction of a regular $3$-gon and $4$-gon using a regular hexagon}
	\label{figure9}
\end{figure}

It is interesting that the construction of a regular hexagon   plays an important role in constructing  other regular $n$-gons. In fact, by connecting odd or even-labeled vertices we always get a regular $3$-gon (equilateral triangle). We can also use a regular hexagon to construct a regular $4$-gon (square) if we use the first and third vertices and the midpoints of the second and fifth arcs of the circumscribed circle (see \cref{figure9}).

\begin{figure}[H]
 	\centering
 	\includegraphics[scale=1]{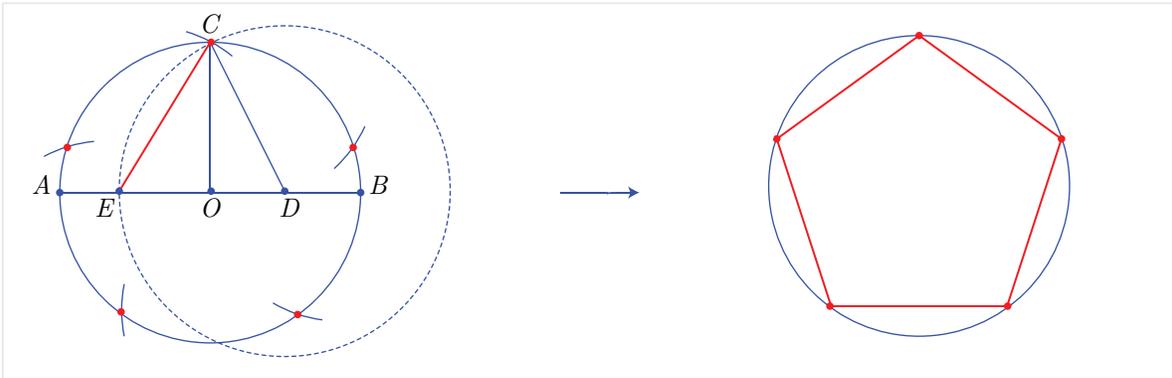}
 	\caption{Ptolemy's construction of a regular pentagon}
 	\label{figure10}
 \end{figure}

A regular pentagon is also    constructible using a compass and straightedge.  However it is not as easy as for a  hexagon.  The main part in this process is to construct the Golden ratio  $\frac{\sqrt{5}-1}{2} $.  Greek mathematicians gave different methods for this problem among which Ptolemy's method seems rather simple. In this method, we start  with a circle of radius $r$ and consider  the diagonal $AB$ and the point $C$ which is the midpoint of the upper arc $ AB$. Then we connect  $C$ to the midpoint of $OB$, say $D$. Name     the intersection point of the circle centered at $D$ with radius $CD$ with radius $AO$ as $E$. It is easy to show that 
\[ \frac{\overline{OE}}{\overline{OB}}= \frac{\overline{EB}}{\overline{OE}}=\frac{\sqrt{5}-1}{2}\]
and the  line segment $CE$  is the side of the inscribed pentagon. By using a compass, we can mark five points on the circle separating five equal arcs and construct a regular pentagon (see \cref{figure10}).

\begin{figure}[H]
	\centering
	\includegraphics[scale=1]{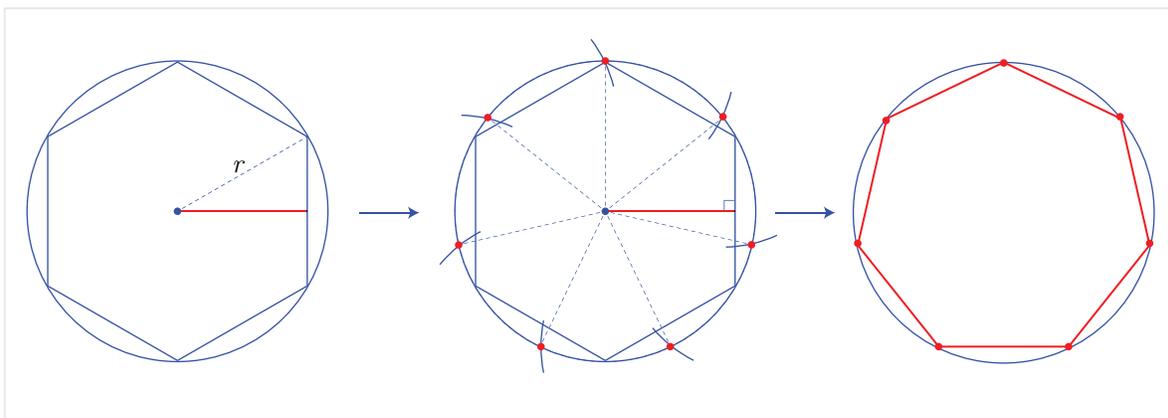}
	\caption{Heron's construction of a regular heptagon}
	\label{figure11}
\end{figure}

So far, we can construct any regular $n$-gon inscribed in a given circle for all \textit{regular numbers} $n= 2^{p}3^{q}5^{r} \geq 3 $, where $p,q,r$ are non-negative integers.  The first number that is not in this list is the prime number 7. That might be the reason that the construction of a regular heptagon has been one of the most interesting problems in the history of mathematics and many mathematicians have tried to tackle it. This problem had obsessed  mathematicians for a long time until the German mathematician Gauss  showed that it is impossible to construct a heptagon with a compass and a straightedge. So, all the constructions before Gauss were approximations, although they were considered   accurate for a long time. As usual, Greek mathematicians such as  Archimedes  and Heron  have played a major role  in  developing  this  theory and many others including Muslim mathematicians followed them to give their own constructions\footnote{See \cite{Hog84}, for a history and list of such constructions.}. Although Archimedes's method is not simple\footnote{See  \cite{Hol10}, for a modern explanation of Archimedes's construction.}, the construction of Heron is easy and practical. He uses a regular hexagon inscribed in a circle of radius $r$ and takes the apothem of  the  regular hexagon as the side of the regular heptagon (see \cref{figure11}). The length of this line segment is $\frac{\sqrt{3}}{2}r$ and the error of the approximation is around $2\%$. Heron uses the approximation $\sqrt{3}=\frac{7}{4} $ to compute the side of the regular heptagon as $ a=\frac{7}{8} r$. 

\begin{figure}[H]
	\centering
	\includegraphics[scale=1]{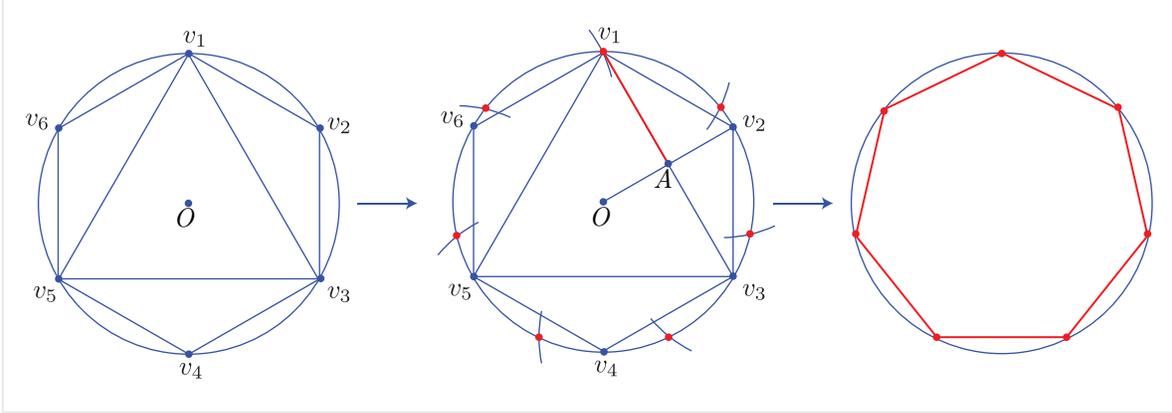}
	\caption{Dürer's construction of a regular heptagon}
	\label{figure12}
\end{figure}

Besides mathematicians,   artists also  tried to solve this problem one of whom was the German artist Albrecht Dürer. In one of his books on geometry, he gave different ways to construct regular polygons. The method for the regular heptagon is given in \cref{figure12}. It is interesting that he has used an inscribed regular hexagon to construct his regular heptagon. He uses the odd-labeled vertices $v_{1}, v_{2}, v_{3} $ of hexagon to construct an equilateral triangle. Then he draws the line passing through the center $O$ and the vertex $v_{2}$ of the hexagon intersecting the side between vertices $v_{1},   v_{3} $ at a point, say $A$. He finally uses the line segment between $A$ and the vertex $v_{1} $ as the side of the required regular heptagon. This method is similar to Heron's construction, because the side of regular heptagon here is also the common height of the six equilateral triangles forming the regular hexagon.

It is clear that the Susa scribe  of this tablet has  used approximate methods to construct the regular polygons on this clay tablet. As we saw, he uses the approximations to  compute the area  of   regular polygons. In the   case of hexagon, fortunately the approximate value of the    side $a$   and the radius  $r$ of   the circumscribed circle are the same. So  he could have   simply taken advantage of this fact and assumed  the radius of the circumscribed circle  to be  the very side of the regular hexagon in order to   construct the figure.  In other words, although he has used an approximation, he has obtained an accurate construction.

\begin{figure}[H]
	\centering
	\includegraphics[scale=1]{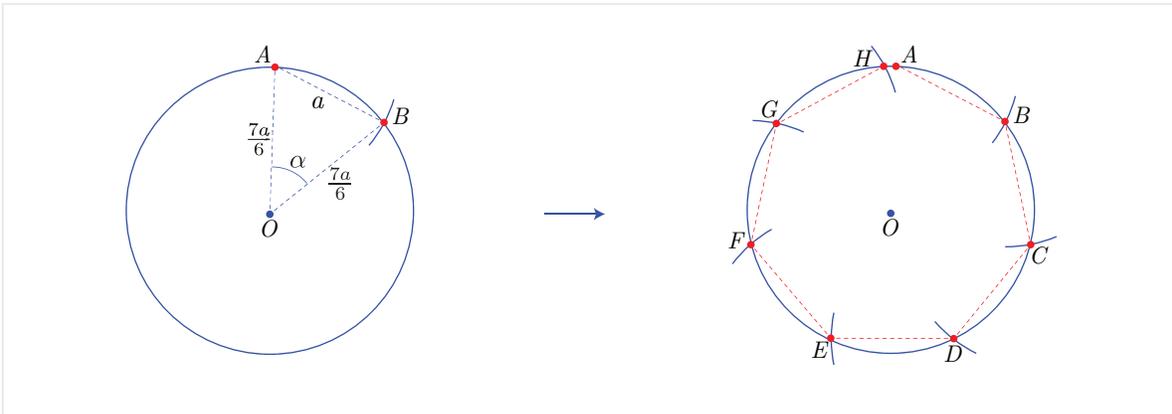}
	\caption{Approximate construction of a regular heptagon}
	\label{figure13}
\end{figure}

In the case of the  heptagon, the  radius of the circumscribed circle   $r$  was  approximated by   $r\approx \frac{7}{6}a$, where $a$ is the side of the regular heptagon. If  he   tried to use this approximate value as the  radius of the circumscribed circle employing  a similar process to that in the case of the  hexagon, as  can be seen in \cref{figure13}, the last arc centered at $G$ would have intersected the circle at  point $H$  not coinciding with the initial point $A$  and there would be a small gap between $H$ and $A$. The cause of this small gap is the approximation $r\approx \frac{7}{6}a$. In fact,   in the isosceles triangle $\triangle OAB$ in \cref{figure13}, we have
$ \sin\left(\frac{\alpha}{2}\right) =\frac{3}{7}$
implying that  the central angle of the  arc  $AB$  is $ \alpha=  2\times \arcsin\left(\frac{3}{7}\right)\approx  50.76^{\circ}$.
So the  total of all central angles   is $7\alpha\approx  355.32^{\circ}$,  which is $4.68^{\circ}$ less than a complete rotation $360^{\circ} $.
 
To complete the figure and get a heptagon, one could  connect  $G$ to either $A$ or   the midpoint of arc $HA $. In both cases, by connecting the adjacent points  of the seven obtained  points, we get an approximation of a regular heptagon (see \cref{figure13}). This might have been the way the Susa scribes   constructed the regular heptagon on the reverse of this tablet.  It should be noted that the   error incurred by using  this method is negligible when one   constructs this figure  on a small clay tablet with dimensions 12cm $\times$ 12cm!

\begin{figure}[H]
	\centering
	\includegraphics[scale=1]{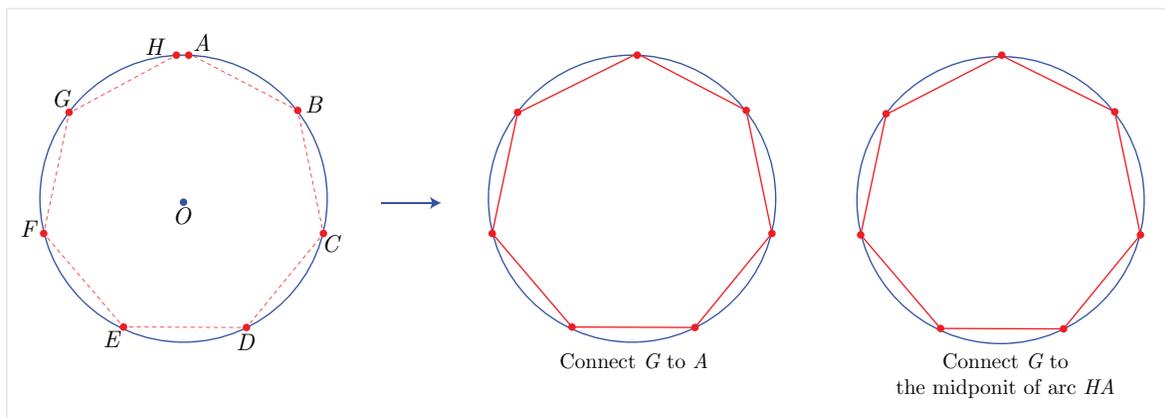}
	\caption{Possible Elamite construction of a regular heptagon}
	\label{figure14}
\end{figure}

\section{Conclusion}
The numerical data  and the figures on    this tablet as well as the related   mathematical calculations    are in our view compelling evidence    that the Susa scribes   applied a standard method to compute the approximate   areas of certain   polygons.  Besides the  standard formulas, they seemed to have used a   simple  and interesting  formula for the area  of a regular heptagon which surprisingly   gives a more accurate approximate value  in comparison to     other  ancient formulas--including  Heron's formula    appearing in the history of mathematics almost 1800 years later. Although the significance of this achievement has not been appreciated  by most scholars, we hope this article sheds  additional  light on the  importance of this formula and enables mathematical historians to reconsider its position in the history of mathematics.
 
The   geometrical figures and their surprisingly accurate dimensions on this tablet suggest  that  Susa scribes  were   familiar with regular polygons and  able to construct them on clay tablets. This characteristic is also acknowledged by other scholars, for example, in \cite{Fri07-1,Fri07-2,Hyp02}. Not only  did  they have the ability to divide a circle into $n$ almost equal sectors for regular numbers, but also  for irregular numbers such as $n=7,11,14$.\footnote{In an ongoing book project on Elamite mathematics, the authors have studied and categorized more than 300 geometric patterns bearing on Elamite artifacts and artworks. Many of these patterns contain specific designs whose structures involve dividing a circle into equal   arcs.} Although they had to use approximate values in irregular cases, it seems that they  had  the  practical skills to carry out  this task as well.


{\small 
\begin{thebibliography}{00000000}
\bibitem[\a'{A}lv20]{Alv20}
{J. \a'{A}lvarez-Mon}, {\em The Art of Elam},  CRC Press,   2020. 

	
\bibitem[BR61]{BR61} 
E. M. Bruins and M. Rutten,    \textit{Textes Math\a'{e}matiques de Suse} [\textbf{TMS}], Librairie Orientaliste Paul Geuthner,   Paris, 1961.

\bibitem[DTKBM]{DTKBM} 
F. Desset, K. Tabibzadeh, , M. Kervran,  G.P.  Basello,   and G. Marchesi,  \textit{The Decipherment of Linear Elamite Writing},   Zeitschrift für Assyriologie und vorderasiatische Archäologie, vol. 112, no. 1, 2022, pp. 11-60. 


	\bibitem[FR98]{FR98}
{D. Fowler  and E. Robson}, {\em Square Root Approximations in Old Babylonian Mathematics: YBC 7289 in Context}, Historia Mathematica  25 (1998), pp. 366-378.


\bibitem[Fri07-1]{Fri07-1}
{J. Friberg}, {\em A Remarkable Collection of Babylonian Mathematical Texts}, Springer, 2007.


\bibitem[Fri07-2]{Fri07-2}
{J. Friberg}, {\em Amazing Traces of a Babylonian Origin in Greek Mathematics}, World Scientific, 2007.

\bibitem[FA16]{FA16}
{J. Friberg and F.  Al-Rawi}, {\em New Mathematical Cuneiform Texts}, Sources and Studies in the History of Mathematics and Physical Sciences, Springer, 2016.

\bibitem[Hea21]{Hea21}
{S. T. Heath}, {\em A History of Greek Mathematics: Volumes I-II}, Oxford University Press,   1921.

\bibitem[Hog84]{Hog84}
{J. P. Hogendijk}, {\em Greek and Arabic constructions of the regular heptagon}, Archive for History of Exact Sciences volume 30 (1984), pp.197--330.

\bibitem[Hol10]{Hol10}
{A. Holme}, {\em Geometry: Our Cultural Heritage}, Springer, Second Edition,   2010.

\bibitem[H\o{}y93]{Hyp93}
{J. H\o{}yrup}, {\em  	Mathematical Susa Texts VII and VIII. A Reinterpretation}, Altorientalische Forschungen 20 (1993), pp. 245--260. 

\bibitem[H\o{}y99]{Hyp99}
{J. H\o{}yrup}, {\em  Pythagorean ``Rule'' and ``Theorem''--Mirror of the Relation Between Babylonian and Greek Mathematics}, published in: Johannes Renger (Ed.), \textit{Babylon: Focus mesopotamischer Geschichte, Wiege früher Gelehrsamkeit, Mythos in der Moderne}
(CDOG 2). Saarbrücker Druckerei und Verlag, Saarbrucken,   (1999), pp. 393--407.

\bibitem[H\o{}y02]{Hyp02}
{J. H\o{}yrup},   {\em  	Lengths, Widths, Surfaces: A Portrait of Old Babylonian  Algebra and Its Kin}, Springer,  Studies and Sources in the History of Mathematics and Physical Sciences, New York,   2002. 

\bibitem[H\o{}y17]{Hyp17}
{J. H\o{}yrup},   {\em  	Algebra in Cuneiform: Introduction to an Old Babylonian  Geometrical Technique}, Max Planck Institute for the History of Science, Textbook 2,  Berlin,   2017. 



	

\bibitem[Mur92-1]{Mur92-1}
{K. Muroi}, {\em  Reexamination of Susa Mathematical Text No. 3: Alleged Value $\pi=3\frac{1}{8}$}, Historia Scientiarum Vol. 2-1 (1992), pp. 45--49.

\bibitem[Mur92-2]{Mur92-2}
{K. Muroi}, {\em  A New Interpretation of Susa Mathematical Text No. 7}, Studies in Babylonian Mathematics No. 2  (1992), pp. 14--19.

\bibitem[Mur92-3]{Mur92-3}
{K. Muroi}, {\em Reexamination of Susa Mathematical Texts Aa (No. 5) and Bb (No. 6)}, Studies in Babylonian Mathematics No. 2  (1992), pp. 22--26.


\bibitem[Mur94-1]{Mur94-1}
{K. Muroi}, {\em  Reexamination of the First Problem of the Susa Mathematical Text No. 9}, Historia Scientiarum Vol. 3-3 (1994), pp. 231--233.

\bibitem[Mur94-2]{Mur94-2}
{K. Muroi}, {\em  Reexamination of   the Susa Mathematical Text No. 8}, Journal of History of Mathematics Vol. 140 (1994), pp. 50--56.



\bibitem[Mur00-1]{Mur00-1}
{K. Muroi}, {\em  Quadratic Equations in   the Susa Mathematical Text No. 21}, SCIAMVS Vol. 1 (2000), pp. 3--10.

\bibitem[Mur01-1]{Mur01-1}
{K. Muroi}, {\em  Inheritance Problems in the Susa Mathematical Text No. 26}, Historia Scientiarum Vol. 10-3 (2001), pp. 226--234.

\bibitem[Mur01-2]{Mur01-2}
{K. Muroi}, {\em  Reexamination of  the Susa Mathematical Text No. 12}, SCIAMVS Vol. 2 (2001), pp. 3--8.



\bibitem[Mur03-2]{Mur03-2}
{K. Muroi}, {\em Excavation Problems  in Babylonian Mathematics Susa Mathematical Text No. 24 and Others}, SCIAMVS Vol. 4 (2003), pp. 3--21.







\bibitem[Mur13]{Mur13}
{K. Muroi}, {\em  Mathematics Hidden Behind the Two Coefficients of Babylonian Mathematics}, 	arXiv:1305.0865 [math.HO], 2013, 10 pages: \url{https://arxiv.org/pdf/1305.0865.pdf}. 



\bibitem[Mur16]{Mur16}
{K. Muroi}, {\em  The Oldest Example of  $\pi=3\frac{1}{8}$  in Sumer: Calculation of the Area of a Circular Plot}, arXiv:1610.03380 [math.HO], 2016, 6 pages: \url{https://arxiv.org/ftp/arxiv/papers/1610/1610.03380.pdf}.




\bibitem[Mur]{Mur}
{K. Muroi}, {\em  Reexamination of  the Susa Mathematical Text No. 11: A Primitive Indeterminate Equation and a Complicated System of Quadratic Equations}, To be published in Acta Sumerologica No. 21, 10 pages.




\bibitem[PAT92]{PAT92}
{H. O. Prudence, J. Aruz, and F. Tallon}, {\em Royal City of Susa: Ancient Near Eastern Treasures in the Louvre}, The Metropolitan Museum of Art, New York, 1992.


 \end{thebibliography} 
 
 }
\end{document}